\documentclass[10pt]{amsart}

\usepackage[T1]{fontenc}
\usepackage{lmodern}
\usepackage{epigraph}
\usepackage{amsmath,amssymb}
\usepackage[numbers]{natbib}
\usepackage{hyperref}

\title{In Memoriam: Professor Dr. Winfried Bruns (1946--2026)}
\author{Peyman Nasehpour}
\email{nasehpour@gmail.com}
\begin{document}

	\begin{abstract}
	Professor Winfried Bruns was a prominent mathematician with vast knowledge in various fields of mathematics. He specialized in combinatorial commutative algebra, algebraic geometry, computer algebra, algebraic coding theory, and cryptography. In his classes, he incorporated the history of mathematics, as well as topics related to music, making his lectures engaging and enriching from an educational perspective. Additionally, he was a skilled musician, specializing in the viola. The main purpose of this note is to explore his life as a mathematician, computer scientist, and musician.
\end{abstract}
	
	\maketitle

	\vspace{1em}
	
	\epigraph{
		Es ist wahr, ein Mathematiker, der nicht etwas Poet ist, wird nimmer ein vollkommener Mathematiker sein.
	}{
		Karl Weierstrass (1815--1897)
	}
	
	\vspace{1.5em}
	
	\section*{Introduction}
	
	Professor Winfried Bruns was, and continues to be, a luminary in the field of commutative algebra. I reaped the benefits of his wisdom and precision. During my doctoral program, I had the privilege of entering his office at any time to ask my questions. Whenever I posed spontaneous questions in commutative algebra, he always provided insightful comments that guided me toward the results I was striving to prove. What I struggled to prove, he foresaw with clarity; it was as if I was at the base of a mountain while he, from the peak, could observe everything clearly. Albert Einstein (1879--1955) once said, ``Der Lehrer hat die Aufgabe, den Sch\"{u}ler zu eigenem Denken und Handeln zu f\"{u}hren,'' which translates to: ``The teacher's task is to lead the student to independent thinking and acting.'' Professor Bruns exemplified such a teacher, based on my own experience. Someone who did not know him personally might think I am exaggerating. However, after checking his biography in this note, they will see that what I have said is, in fact, an understatement.
	
	One pedagogical characteristic of Professor Bruns in teaching and leading students in mathematics was his constant encouragement to pursue the subject. Once, I discussed a mathematical problem with him that seemed simple but was actually quite difficult. He encouraged me to tackle the problem. I expressed my doubts, saying, ``If knowledgeable and experienced mathematicians haven't been able to fully solve it, how can I?'' His response surprised me: ``Sometimes, a fresh mind may find the solution easier than others.''
	
	Professor Bruns was exceptionally caring towards his students. Whenever we met at the Institute for Mathematics and Informatics at the University of Osnabr\"{u}ck, he always inquired about our well-being. On one occasion, when a storm was forecasted in Osnabr\"{u}ck, he personally called us to ensure we went home early and didn't stay at the Institute. Additionally, when Matteo Varbaro joined our Institute as a guest researcher for a month, I observed that Professor Bruns was consistently attentive and took great care of him.
	
	Professor Bruns was a distinguished expert in computer algebra and applied mathematics as well. Not only was he the main person behind the software Normaliz, but his lecture notes on algebraic coding theory and cryptography also demonstrate his profound knowledge in applied algebra. Together with his former supervisor Professor Udo Vetter, he authored a two-volume book on calculus in German, published by the University of Oldenburg, showcasing their talent in mathematics education \cite{BrunsVetter2003}.
	
	Professor Bruns' knowledge of the history of mathematics, particularly algebra, was remarkable. Some years ago, I inquired over the phone about the history of the phrase ``determinant trick'' in commutative algebra. He provided a hint that helped me trace the use and history of this phrase across various resources. On another occasion, I asked his opinion about a quote by a mathematician. He referred me to a book on the history of mathematics that I had never heard of, cautioning that one should carefully examine the claims made in such books. Another valuable lesson I learned from him was during a class at the University of Osnabr\"{u}ck in 2008. He explained that while many parts of the extensive series on mathematics by Nicolas Bourbaki, originally written in French, have been translated into English, not all chapters had been translated at that time. The situation may have changed since then.
	
	Once, in the cafeteria of the Department of Cognitive Science at the University of Osnabr\"{u}ck, I had a discussion about the pronunciation of some German words. During this conversation, I discovered that he had considerable knowledge in phonetics. At that time, I was surprised, but later I learned that he had attended the classical gymnasium in Goslar, which focused on ancient languages. Perhaps his interest in different languages stemmed from those years. Being fluent in German and English and having a good knowledge of French and Italian, the books he wrote on mathematics in German demonstrate his deep understanding of German literature. One of his favorite texts was ``Die Blechtrommel'' (The Tin Drum), a novel by G\"{u}nter Grass, first published in 1959. Additionally, in his report on my dissertation, which I wrote under his supervision, he used the German noun ``Findigkeit.'' Interestingly, I later discovered that many graduate students, whose mother tongue was German, had never heard the word and had no idea what it meant. The German adjective ``findig'' describes someone who is clever and able to find solutions to problems.

	In addition to his academic achievements, Professor Bruns was an accomplished viola player and performed with the symphonic orchestra of Osnabr\"{u}ck University for many years. On one occasion, he kindly gave me a ticket to the orchestra's performance in a grand church in Osnabr\"{u}ck. As a professional musician, I was deeply impressed by the quality of their performance. Some years earlier, I attended another concert by the university's symphonic orchestra at the Aula in the university's Schloss. Professor Bruns played the viola, and his daughter Julia was the first violinist. The concert was truly high-class, well-performed, and enjoyable even for me with my oriental ears. This is an aspect of his character that was rarely mentioned. Last but not least, it is not surprising to see that Professor Bruns also gave lectures on the relationship between music and mathematics.
	
	\section*{A Brief Biography of Professor Winfried Bruns}
	
	Born on May 5, 1946, in Oker near Goslar, a historic town in Lower Saxony, Germany, Winfried Bruns attended the classical gymnasium in Goslar, which focused on ancient languages---Latin and Greek---from 1956 to 1965. During this time, he also studied natural sciences, including physics and mathematics. After graduating, he served in the army for two years. In 1967, he entered TH Hannover (today, TU Hannover) to study physics, as he had initially dreamed of becoming a physicist. However, he soon developed a greater interest in mathematics and decided to continue his studies in this field. Throughout his studies and until he received his diploma, he was supported by a distinguished grant for excellent students from the Studienstiftung des Deutschen Volkes.
	
	After graduating, he became a scientific assistant at TU Clausthal, where he wrote his doctoral dissertation, Beispiele reflexiver Differentialmoduln \cite{Bruns1972}, under the supervision of Prof. Dr. Udo Vetter (1938--2023). Large parts of his dissertation were later published in the ``Journal f\"{u}r die reine und angewandte Mathematik'' under the title ``Zur Reflexivit\"{a}t analytischer Differentialmoduln'' in 1975 \cite{Bruns1975}. At that time, ``Basic Element Theory'' had gained popularity as an important tool for studying generating sets of modules. By extending techniques on basic elements, as developed by Eisenbud and Evans in their influential paper ``Basic elements: Theorems from Algebraic $K$-Theory,''(cf. \cite{EisenbudEvans1972}) Bruns proved a beautiful and surprising theorem in his paper ``Jede'' endliche freie Aufl\"{o}sung ist freie Aufl\"{o}sung eines von drei Elementen erzeugten Ideals \cite{Bruns1976}. This marked the beginning of an outstanding scientific career. The methods used in this paper, specifically those of basic element theory, were also the subject of his habilitation dissertation in 1977. Notably, the project titled ``Factoring Out a General Element from an $s$-th Syzygy'' in Eisenbud's celebrated book on commutative algebra \cite{Eisenbud1995} explains that the phenomena most likely to be discovered by the reader were first noticed and exploited by W. Bruns.
	
	The Normalschule in Vechta was established in 1830; the term ``Normalschule'' historically refers to a type of teacher training school in Germany. In 1969, the school in Vechta, along with seven other institutions, united to form the Lower Saxony University of Education. In 1973/74, the Vechta University of Education became a department of the newly founded University of Osnabr\"{u}ck and remained so until 1994. In 1979, Bruns became a full professor at the University of Osnabr\"{u}ck, Division Vechta. The University of Vechta gained independence from the University of Osnabr\"{u}ck and became a separate institution in 1995. However, Bruns continued his professorship at the University of Osnabr\"{u}ck until he retired in 2014. Even after his official retirement, Professor Bruns continued to be active as a Professor Emeritus. For example, in 2022, Winfried Bruns, Aldo Conca, Claudiu Raicu, and Matteo Varbaro authored the book ``Determinants, Gr\"{o}bner Bases and Cohomology,'' offering an up-to-date, comprehensive account of determinantal rings and varieties \cite{BCRV2022}. This book presents a multitude of methods used in their study, incorporating tools from combinatorics, algebra, representation theory, and geometry.
	
	Before his fruitful career in Vechta, he spent a year as a Visiting Lecturer at the University of Illinois at Urbana-Champaign, primarily to collaborate with E. Graham Evans (1942--2021). This cooperation had a formative effect on his work and resulted in a lifelong friendship with Evans. One notable achievement from his visit was the joint publication with Evans and Griffith in 1980 of the paper Syzygies, Ideals of Height Two, and Vector Bundles \cite{BrunsEvansGriffith1980}. This paper addressed the syzygy problem in some special cases, providing affirmative answers. A year later, Evans and Griffith proved the syzygy theorem in general for any Noetherian local domain containing a field \cite{EvansGriffith1981}. In subsequent years, Bruns and others further generalized this theorem in various directions.
	
	During his tenure in Vechta, Prof. Bruns published numerous influential papers on key aspects of modern commutative algebra. His work covered topics such as generic maps, generic resolutions, divisors on varieties of complexes, straightening laws on modules, and straightening closed ideals. The theory of algebras with straightening laws (ASLs), initially developed by De Concini, Eisenbud, and Procesi in the early 1980s, included determinantal rings as prominent examples. Bruns expanded this theory by demonstrating that it provides upper bounds for the arithmetical rank of ideals generated by a poset ideal of the poset underlying the algebra with straightening law. He also introduced modules with straightening laws to study the symmetric algebra of generic modules.
	
	In 1988, Bruns and Vetter published their influential Springer lecture notes titled Determinantal Rings, which became an enchiridion on the subject \cite{BrunsVetter1988}. Remaining the standard reference in the field, the book covers a wide range of topics, including the classical theory developed by Hochster and Eagon and the powers and symbolic powers of determinantal ideals. Additionally, the theory of algebras with straightening laws (ASLs) is employed and the canonical class of determinantal rings is determined, among many other significant results.
	
	Gr\"{o}bner bases entered the theory of determinantal ideals with the work of Sturmfels in 1990. This aspect of the theory had a significant impact on Bruns' later work and his future scientific collaborations. In 2003, jointly with Aldo Conca, he wrote the paper Gr\"{o}bner Bases and Determinantal Ideals \cite{BrunsConca2003}, continuing the line of investigation started by Sturmfels. Following their first joint paper, KRS (Knuth-Robinson-Schensted) and powers of determinantal ideals in 1998, Bruns and Conca enjoyed an extremely fruitful collaboration. They authored a series of influential papers on the theory of determinantal rings, studying Gr\"{o}bner bases and powers of determinantal ideals. In other papers, jointly with Aldo Conca and Tim R\"{o}mer, they explored Koszul cycles, Koszul homology, and syzygies of Veronese subalgebras \cite{BrunsConcaRoemer2011}. Later, Matteo Varbaro joined the team, and in 2013, Bruns, together with Conca and Varbaro, published a fundamental paper in Advances in Mathematics addressing the complex problem of understanding the relations of the algebra generated by the $t$-minors of a matrix of indeterminates \cite{BrunsConcaVarbaro2013}.
	
	Bruns' book Cohen-Macaulay Rings\cite{BrunsHerzog1993} (coauthored with J\"{u}rgen Herzog) became highly popular and is considered one of the primary resources in commutative algebra. Alongside David Eisenbud's Commutative Algebra with a View Toward Algebraic Geometry, this monograph remains essential reading for any advanced student in the field.
	
	His other book, Polytopes, Rings, and $K$-theory \cite{BrunsGubeladze2009}, coauthored with Joseph Gubeladze, explores the interplay between discrete convex geometry, commutative ring theory, algebraic $K$-theory, and algebraic geometry, showcasing his extensive knowledge across various fields of mathematics. This book also marked the culmination of years of scientific collaboration with Gubeladze, which began in 1996.
	
	In one of their notable papers on the normality and covering properties of affine semigroups, published in 1999 in the Journal f\"{u}r die reine und angewandte Mathematik, they discovered a six-dimensional counterexample to a conjecture by Seb\"{o} \cite{BrunsGubeladze1999}. Seb\"{o} had conjectured that a finite rational cone admits a unimodular covering by simplicial cones spanned by elements of the Hilbert basis. It was previously known by Seb\"{o} and others that three-dimensional rational cones could even admit unimodular triangulations, but such triangulations generally do not exist in dimension $4$. In their paper, Bruns and Gubeladze also presented an algorithm to determine whether a finite rational cone admits a unimodular covering. The computational and algorithmic aspects of the theory of polytopes and affine monoids were always of particular interest to Bruns.
	
	In the same year and journal, Bruns, along with Gubeladze, Henk, Martin, and Weismantel, demonstrated that the Bruns--Gubeladze counterexample to unimodular coverings also served as a counterexample to the integral Carathéodory property of cones \cite{BGHMW1999}. This property requires that any integral vector of an n-dimensional integral polyhedral pointed cone can be expressed as a nonnegative integral combination of at most n elements of the Hilbert basis of the cone. Their result showed that Carath\'{e}odory's theorem for convex cones does not have an integer analogue.
	
	The discovery of these counterexamples significantly advanced experimental mathematics. This progress has been facilitated by the immense capacity of modern computers and powerful computer algebra systems, such as Normaliz. Normaliz, created by Bruns and his former doctoral student Robert Koch, was regularly updated by him and his team, which has included Ichim, Sieg, R\"{o}mer, and S\"{o}ger. The program is freely downloadable and integrated into major algebraic software suites such as Macaulay2, polymake, SageMath, and Singular. Up to now, it remains an indispensable tool in the study of polytopes and affine monoids.
	
	As a primary supervisor, Professor Bruns mentored eight doctoral students: Margherita Barile, Stefan G\"{u}nther, Robert Koch, Peyman Nasehpour, Richard Sieg, Christof S\"{o}ger, Jan Uliczka, and Michael von Thaden. Additionally, he served as a secondary supervisor for Bogdan Nicolae Ichim and Klaus G. Warneke. Beyond these formal roles, he supported numerous young students, guiding them into the field of mathematics. In recognition of his dedication to teaching, he was awarded the Osnabr\"{u}ck University Prize for Excellence in Teaching in 2002.
	
	Together with Holger Brenner, Tim R\"{o}mer, and other colleagues, Bruns successfully applied for a Graduiertenkolleg from the DFG (German Research Foundation). Established in 2013, this Research Training Group has been highly successful, providing doctoral and postdoctoral positions for many young researchers. His efforts significantly contributed to making the mathematics department at Osnabr\"{u}ck one of the leading centers for commutative algebra in Europe.
	
	Over the course of his career, Prof. Winfried Bruns published at least 126 papers, including a number of preprints on arXiv. He also served on the editorial boards of Journal of Communications in Algebra, and Journal of Homology, Homotopy and Applications. His work profoundly shaped and permanently influenced the fields of commutative algebra and computational mathematics.
	
	On May 30th, I called him and we spoke. He was delighted that his enduring contributions to mathematics were to be celebrated by colleagues and students from around the world at a dedicated conference in Osnabr\"{u}ck, scheduled for June 24--26, 2026. Tragically and unexpectedly, just days before the community could gather to honor his monumental achievements in mathematics, Professor Bruns passed away in Genoa, Italy, on 21 June, 2026, following a sudden stroke. The conference, originally conceived to celebrate Winfried Bruns' 80th birthday, was ultimately held in his memory.
	
	\section*{Acknowledgments}
	
	Some parts of the biography of Prof. Dr. Winfried Bruns are based on information originally provided by Prof. Dr. J\"{u}rgen Herzog (1941--2024) (see \cite{CGR2017}).

	\end{document}